%%ineqcs_arxiv.tex 
\input amstex 
\documentstyle{amsppt} 
%\nologo
\loadbold

\magnification=1200
\hsize=5.75truein
\vsize=8.75truein 
\hcorrection{.25truein}
\loadeusm \let\scr\eusm

\font\Rm=cmr12

\define\Aut#1#2{\text{\rm Aut}_{#1}(#2)}
\define\End#1#2{\text{\rm End}_{#1}(#2)}
\define\Hom#1#2#3{\text{\rm Hom}_{#1}({#2},{#3})} 
\define\GL#1#2{\roman{GL}_{#1}(#2)}
\define\M#1#2{\roman M_{#1}(#2)}
\define\Gal#1#2{\text{\rm Gal\hskip.5pt}(#1/#2)}
\define\upr#1#2{{}^{#1\!}{#2}} 
\define\pre#1#2{{}_{#1\!}{#2}} 
\define\sw{\text{\rm Sw}} 
\define\ar{\text{\rm Ar}}
\define\sve{{\ssize\vee\,}}
\let\ge\geqslant
\let\le\leqslant

\let\ve\varepsilon

\let\vS\varSigma 
\let\vs\varsigma

\define\wW#1{\widehat{\scr W}^{\text{\rm irr}}_{#1}} 
\define\wss#1{\widehat{\scr W}^{\text{\rm ss}}_{#1}} 
\define\dw#1{\widehat{\scr W}^\scr D_{#1}} 
\font\sns=cmss10 \define\Tau{{\text{\sns T}}} 
\define\Rho{{\text{\sns R}}} 
%%%%%%%%%%%%%%%%%%%%%%%%%%%%%%%%%%%%%%%%%%%%%%%%%%
%\document \baselineskip=14pt \parskip=4pt plus 1pt minus 1pt 
\topmatter \nologo \nopagenumbers
\title 
Strong exponent bounds for the local Rankin-Selberg convolution 
\endtitle 
\rightheadtext{Bounds for Rankin-Selberg exponents} 
\author 
Colin J. Bushnell and Guy Henniart 
\endauthor 
\leftheadtext{C.J. Bushnell and G. Henniart}
\affil 
King's College London and Universit\'e de Paris-Sud 
\endaffil 
\address 
King's College London, Department of Mathematics, Strand, London WC2R 2LS, UK. 
\endaddress
\email 
colin.bushnell\@kcl.ac.uk 
\endemail
\address 
Laboratoire de Math\'ematiques d'Orsay, Univ\. Paris-Sud, CNRS, Universit\'e
Paris-Saclay, 91405 Orsay, France.
\endaddress 
\email 
Guy.Henniart\@math.u-psud.fr 
\endemail 
\date March 2016 \enddate 
\abstract 
Let $F$ be a non-Archimedean locally compact field. Let $\sigma$ and $\tau$ be finite-dimensional semisimple representations of the Weil-Deligne group of $F$. We give strong upper and lower bounds for the Artin and Swan exponents of $\sigma\otimes\tau$ in terms of those of $\sigma$ and $\tau$. We give a different lower bound in terms of $\sigma\otimes\check\sigma$ and $\tau\otimes\check\tau$. Using the Langlands correspondence, we obtain the bounds for Rankin-Selberg exponents. 
\endabstract 
\thanks 
The first-named author thanks King's College London for support during the preparation of this paper. 
\endthanks 
\endtopmatter 
\document \baselineskip=14pt \parskip=4pt plus 1pt minus 1pt
%%%%%%%%%%%%%%%%%%%%%%%%%%%%%%%%%%%%%%%%%%%%%%%%%% 
\subhead 
1 
\endsubhead 
Let $F$ be a non-Archimedean, locally compact field. For integers $m,n\ge 1$ let $\pi$, $\rho$ be irreducible, smooth, complex representations of the general linear groups $\GL mF$, $\GL nF$ respectively. If $s$ is a complex variable and $\psi$ a non-trivial smooth character of $F$, we consider the $\roman L$-function $\roman L(\pi\times\rho,s)$ and the local constant $\ve(\pi\times\rho,s,\psi)$ of \cite{17} or \cite{21}, \cite{22}. If $q$ is the cardinality of the residue field of $F$, the local constant takes the form 
$$ 
\ve(\pi\times\rho,s,\psi) = \ve(\pi\times\rho,0,\psi)\,q^{-s(\ar(\pi\times\rho)+mnc(\psi))}. 
$$ 
Here, $c(\psi)$ is an integer depending only on $\psi$. The integer $\ar(\pi\times\rho)$ depends only on the pair $(\pi,\rho)$. Here we call it the {\it Rankin-Selberg exponent\/} of $(\pi,\rho)$. 
\par 
If we take $n=1$ and let $\rho$ be the trivial character $1$ of $F^\times \cong \GL1F$, then $\ve(\pi\times1,s,\psi)$ is the Godement-Jacquet local constant $\ve(\pi,s,\psi)$ \cite{10}, and $\ar(\pi\times1)$ is denoted simply $\ar(\pi)$. The aim of this paper is to give strong, universal estimates for $\ar(\pi\times\rho)$ in terms of $\ar(\pi)$ and $\ar(\rho)$. We give a second lower bound in terms of exponents of the pairs $(\pi,\check\pi)$, $(\rho,\check\rho)$. These results are Corollaries A--C below.  
\subhead 
2 
\endsubhead 
We fix a separable algebraic closure $\bar F/F$ of the field $F$, and form the Weil group $\scr W_F = \scr W_{\bar F/F}$. Let $\dw F$ be the set of equivalence classes of finite-dimensional, semisimple representations of the Weil-Deligne group defined by $\scr W_F$. With $\pi$ and $\rho$ as before, the Langlands correspondence \cite{11}, \cite{15}, \cite{18}, \cite{19} associates to $\pi$, $\rho$ representations $\upr L\pi, \upr L\rho \in \dw F$. These have dimension $m$, $n$ respectively. 
\par 
For $\sigma\in \dw F$, of dimension $d$, let $\ve(\sigma,s,\psi)$ be the Langlands-Deligne local constant \cite{3},\cite{23} of $\sigma$. Again, $\ve(\sigma,s,\psi) = \ve(\sigma,0,\psi)\,q^{-s(\ar(\sigma)+dc(\psi))}$ and the integer $\ar(\sigma)$ is the {\it Artin exponent\/} of $\sigma$. A defining property of the Langlands correspondence \cite{14}, \cite{16} is that 
$$ 
\ve(\pi\times\rho,s,\psi) = \ve(\upr L\pi\otimes \upr L\rho,s,\psi). 
$$ 
Consequently, $\ar(\pi\times\rho) = \ar(\upr L\pi\otimes\upr L\rho)$ and $\ar(\pi) = \ar(\upr L\pi)$. We may therefore tackle the Rankin-Selberg exponent via the Artin exponent of tensor products of representations of the Weil-Deligne group. 
\subhead 
3 
\endsubhead 
We state our results for representations of the Weil-Deligne group. If $\sigma\in \dw F$, $\sigma\neq 0$, write 
$$ 
\eta(\sigma) = \ar(\sigma)/\dim\sigma, \quad \sigma \in \dw F,\ \sigma\neq0. 
$$ 
\remark{Convention} 
When $\sigma$ is the zero representation, $\eta(\sigma)$ is undefined. So, use of the symbol $\eta(\sigma)$ here will {\it always\/} entail the implicit assumption $\sigma\neq 0$. \
\endremark 
Say that $\sigma$ is {\it $\eta$-minimal\/} if $\ar(\sigma) \le \ar(\chi\otimes\sigma)$, for any character $\chi$ of $\scr W_F$. 
\proclaim{Theorem A} 
If $\sigma\in \dw F$ is $\eta$-minimal, then 
$$ 
\eta(\sigma\otimes\tau) \ge \tfrac12\,\roman{max}\,\{\eta(\sigma),\eta(\tau)\},  
$$ 
for all $\tau \in \dw F$. 
\endproclaim 
A trivial example shows that some hypothesis of minimality is required for a result of this kind: for fixed $\sigma,\tau\in \dw F$ and a character $\chi$ of $\scr W_F$, one has $\eta((\chi\otimes\sigma)\otimes (\chi^{-1}\otimes\tau)) = \eta(\sigma\otimes\tau)$. For suitable choice of $\chi$, one has $\eta(\chi\otimes\sigma) = \eta(\chi^{-1}\otimes\tau) = \eta(\chi)$ and this may be taken as large as desired. 
\par 
Further examples show that the constant $\frac12$ is best possible: there are many pairs of {\it irreducible\/} representations $(\sigma, \tau)$, with $\sigma$ being $\eta$-minimal, for which $2\eta(\sigma\otimes\tau) = \eta(\sigma) = \eta(\tau)$. However, by restricting the class of representation one can get better constants: see the examples in 2.3 and 3.5. 
\par 
There is a second, rather different, lower bound. This avoids the necessity for a minimality condition by using the operation $\sigma\mapsto \check\sigma$ of contragredience on $\dw F$. 
\proclaim{Theorem B} 
If $\sigma,\tau\in \dw F$, then 
$$ 
\eta(\sigma\otimes\check\tau) \ge \tfrac12\big(\eta(\sigma\otimes\check\sigma) + \eta(\tau\otimes\check\tau)\big). 
$$ 
If $\sigma$ and $\tau$ are indecomposable, then $\eta(\sigma\otimes\check\tau) \ge \roman{max}\,\{\eta(\sigma\otimes\check\sigma), \eta(\tau\otimes\check\tau)\}$. 
\endproclaim 
The easy example $\sigma = \tau$ shows that the constant $\frac12$ is again best possible. With regard to upper bounds, we prove: 
\proclaim{Theorem C} 
Let $\sigma,\tau\in \dw F$ have dimensions $m$, $n$ respectively. The Artin exponent $\ar(\sigma\otimes\tau)$ satisfies 
$$ 
\ar(\sigma\otimes\tau) \le n\ar(\sigma) + m\ar(\tau) - \roman{min}\,\{\ar(\sigma),\ar(\tau)\}. 
$$ 
If $\sigma,\tau\in \dw F$ are irreducible, then 
$$ 
\eta(\sigma\otimes\tau) \le \roman{max}\,\{\eta(\sigma),\eta(\tau)\}. 
$$
\endproclaim 
Both aspects of the result are best possible. 
\subhead 
4 
\endsubhead 
Let $\pi$, $\rho$ be irreducible, smooth, complex representations of $\GL mF$, $\GL nF$ respectively. Set $\eta(\pi\times\rho) = \ar(\pi\times\rho)/mn$ and $\eta(\pi) = \ar(\pi)/m$, with the same convention regarding zero representations. Say that $\pi$ is {\it $\eta$-minimal\/} if $\eta(\pi) \le \eta(\chi\pi)$ for all characters $\chi$ of $F^\times$. The Langlands correspondence respects contragredience and twisting with characters, so we have the following consequences of Theorems A--C. 
\proclaim{Corollary A} 
Let $\pi$, $\rho$ be irreducible representations of  the groups $\GL mF$, $\GL nF$ respectively. If $\pi$ is $\eta$-minimal, then 
$$ 
\eta(\pi\times\rho) \ge \tfrac12\,\roman{max}\,\{\eta(\pi),\eta(\rho)\}. 
$$ 
\endproclaim 
\proclaim{Corollary B} 
If $\pi$, $\rho$ are irreducible representations of the groups $\GL mF$, $\GL nF$ respectively, then 
$$ 
\eta(\pi\times\check\rho) \ge \tfrac12\,\big(\eta(\pi\times\check\pi)+\eta(\rho\times\check\rho)\big). 
$$ 
If $\pi$ and $\rho$ are essentially square-integrable, then 
$$ 
\eta(\pi\times\check\rho) \ge \roman{max}\,\{\eta(\pi\times\check\pi), \eta(\rho\times \check\rho)\}. 
$$ 
\endproclaim 
\proclaim{Corollary C} 
Let $\pi$, $\rho$ be irreducible representations of the groups $\GL mF$, $\GL nF$ respectively. The Rankin-Selberg exponent satisfies 
$$ 
\ar(\pi\times\rho) \le n\ar(\pi)+m\ar(\rho) - \roman{min}\,\{\ar(\pi),\ar(\rho)\}. 
$$ 
If the representations $\pi$ and $\rho$ are cuspidal, then 
$$
\eta(\pi\times\rho) \le \roman{max}\,\eta(\pi),\eta(\rho)\}. 
$$ 
\endproclaim 
In all of these statements, the representations $\pi$, $\rho$ are assumed {\it smooth.} Corollary C may also be found in \cite{1}, where it receives a different proof. 
\par 
Beyond remarking that the representation $\pi$ is essentially square-integrable (resp\. cuspidal) if and only if $\upr L\pi \in \dw F$ is indecomposable (resp\. irreducible), there is nothing more to be said about these corollaries. 
\subhead 
5 
\endsubhead 
We return to the Galois side. Let $\wss F$ be the set of equivalence classes of finite-dimensional, smooth, semisimple representations of $\scr W_F$. 
\par 
There is a parallel, but distinct, family of estimates governing the {\it Swan exponent\/} $\sw(\sigma)$, $\sigma\in \dw F$, in place of the Artin exponent. We include them here since, in applications, the Swan exponent often occurs more naturally than the Artin exponent and it can be bothersome to switch between the two languages. The exponent $\sw(\sigma)$ depends only on the restriction of $\sigma$ to $\scr W_F$, so nothing is lost by treating $\sw$ as a function on $\wss F$. 
\par 
If $\sigma \neq 0$, we set $\vs(\sigma) = \sw(\sigma)/\dim\sigma$. Again, use of the symbol $\vs(\sigma)$ entails the implicit assumption $\sigma\neq 0$. 
\par 
Say that $\sigma$ is {\it $\vs$-minimal\/} if $\vs(\sigma) \le \vs(\chi\otimes\sigma)$ for all characters $\chi$ of $\scr W_F$. (Note that the concepts of $\eta$-minimality and $\vs$-minimality are distinct.) We then have the following results. 
\proclaim{Theorem AS} 
If $\sigma\in \wss F$ is $\vs$-minimal, then 
$$ 
\vs(\sigma\otimes\tau) \ge \tfrac12\,\roman{max}\,\{\vs(\sigma),\vs(\tau)\}, 
$$ 
for all $\tau \in \wss F$. 
\endproclaim 
\proclaim{Theorem BS} 
If $\sigma,\tau\in \wss F$, then 
$$ 
\vs(\sigma\otimes\check\tau) \ge \tfrac12\big(\vs(\sigma\otimes\check\sigma) + \vs(\tau\otimes\check\tau)\big). 
$$ 
If $\sigma$ and $\tau$ are irreducible, then $\vs(\sigma\otimes\check\tau) \ge \roman{max}\,\{\vs(\sigma\otimes\check\sigma), \vs(\tau\otimes\check\tau)\}$. 
\endproclaim 
\proclaim{Theorem CS} 
Let $\sigma,\tau\in \wss F$ have dimensions $m$, $n$ respectively. The Swan exponent $\sw(\sigma\otimes\tau)$ satisfies 
$$ 
\sw(\sigma\otimes\tau) \le n\,\sw(\sigma) + m\,\sw(\tau) - \roman{min}\,\{\sw(\sigma),\sw(\tau)\}. 
$$ 
If $\sigma,\tau\in \dw F$ are irreducible, then $\vs(\sigma\otimes\tau) \le \roman{max}\,\{\vs(\sigma),\vs(\tau)\}$. 
\endproclaim 
\subhead 
6 
\endsubhead 
We review some background material in section 1. The proof of Theorem A starts in section 2, where we deal with irreducible representations.  At present, these can only be treated via parallel properties for irreducible {\it cuspidal\/} representations of general linear groups and then using the Langlands correspondence. The method relies on the explicit formula for $\ar(\pi\times\rho)$ in \cite{6}, combining the classification theory of \cite{7}, \cite{8}, \cite{9} with the interpretation \cite{22} of the Rankin-Selberg exponent as a relative Plancherel measure. This is where the factor $\frac12$ of Theorem A first appears and reveals itself as best possible. The main part of the proofs of Theorems A and AS is in section 3. The arguments are all conducted on the Galois side. They are essentially elementary although, in places, they feel intricate. 
\par 
Theorem B and BS are treated in section 4. The proofs start from relatively simple properties of tensor products of irreducible representations observed in \cite{12}, \cite{5} but are equally intricate. For the pairs A/AS, B/BS of parallel theorems, the proofs start together. We then concentrate on the more involved case of the Artin exponent. That done, the argument for the Swan exponent follows a shorter version of the same route, obtained by a simple change of vocabulary. We indicate the process briefly at the ends of the relevant sections. The results are not so easy to deduce from each other, and nothing seems to be gained from constructing an artifical framework in which they can be treated together. The proofs of Theorems C and CS are short, and combined in section 5. 
\remark{Acknowledgement} 
We thank Farrell Brumley and Erez Lapid for interesting questions which drew our attention to this area. 
\endremark 
%%%%%%%%%%%%%%%%%%%%%%%%%%%%%%%%%% 
\head\Rm 
1. Representations of the Weil-Deligne group 
\endhead 
We retain the notations $\scr W_F$, $\wss F$ and $\dw F$ of the introduction. Let $\wW F$ be the set of isomorphism classes of {\it irreducible\/} smooth representations of $\scr W_F$. Starting from the discussions in \cite{3} and \cite{23}, we recall some basic features of representations $\sigma\in \dw F$. We define the Artin exponent in terms of the Langlands-Deligne local constant and collect a number of facts and simple results for use in later sections. 
\subhead 
1.1 
\endsubhead 
Let $q$ be the cardinality of the residue class field of $F$. Let $x\mapsto \|x\|$ denote the unique character of $\scr W_F$ that is trivial on the inertia subgroup of $\scr W_F$ and takes the value $q^{-1}$ on geometric Frobenius elements. 
\par 
For our purposes, a representation $\sigma$ of the Weil-Deligne group of $F$ is a pair $(\sigma_\scr W,\frak n)$ consisting of a finite-dimensional, smooth, semisimple representation $\sigma_\scr W:\scr W_F \to \Aut{\Bbb C}V$ and a nilpotent endomorphism $\frak n$ of the vector space $V$ such that 
$$ 
\sigma_\scr W(g)\,\frak n = \|g\|\,\frak n\,\sigma_\scr W(g), \quad g\in \scr W_F. 
$$ 
We denote by $\dw F$ the set of isomorphism classes of such representations. For $\sigma\in \dw F$, we rarely use the notation $\sigma_\scr W$ but speak instead of the restriction of $\sigma$ to $\scr W_F$. In the same spirit, a representation $\sigma\in \wss F$ defines an element $(\sigma,0)$ of $\dw F$ that we continue to denote by $\sigma$. 
\par
The set $\dw F$ admits a notion of direct sum, 
$$ 
(\sigma,\frak n)\oplus (\tau,\frak n) = (\sigma\oplus\tau,\frak m\oplus \frak n). 
$$ 
We say $(\sigma,\frak n)$ is {\it indecomposable\/} if it cannot be expressed in this way as a direct sum in which both factors are non-trivial. Surely any $\sigma\in \dw F$ may be expressed as a direct sum of indecomposable elements of $\dw F$. Such a decomposition is unique up to permutation of the isomorphism classes of indecomposable factors. 
\par 
To define the tensor product $(\sigma,\frak m)\otimes (\tau,\frak n)$, let $\sigma$ act on a vector space $V$ and $\tau$ on $W$. One sets 
$$ 
(\sigma,\frak m)\otimes (\tau,\frak n) = (\sigma\otimes\tau, \frak m\otimes 1_W + 1_V\otimes\frak n). 
$$ 
\subhead 
1.2 
\endsubhead 
We recall the standard first example of an element of $\dw F$. Let $n\ge 1$ be an integer and let $\roman{sp}_n \in \wss F$ denote the direct sum of the characters $x\mapsto \|x\|^i$, for $0\le i\le n{-}1$. We view $\roman{sp}_n$ as acting on $V = \Bbb C^n$. The space $V$ admits a {\it regular\/} nilpotent endomorphism $\frak n$ such that $\roman{Sp}_n(1) = (\roman{sp}_n,\frak n)$ is a representation of the Weil-Deligne group. The isomorphism class of $\roman{Sp}_n(1)$ is independent of the choice of $\frak n$. 
\par 
More generally, let $\sigma\in \wW F$. We define 
$$ 
\roman{Sp}_n(\sigma) = \sigma\otimes\roman{Sp}_n(1). 
$$ 
An exercise \cite{23} (4.1.5) yields: 
\proclaim{Fact} 
A representation $\vS\in \dw F$ is indecomposable if and only if $\vS = \roman{Sp}_n(\sigma)$, for an integer $n\ge 1$ and a representation $\sigma\in \wW F$. Moreover, $\roman{Sp}_n(\sigma) \cong \roman{Sp}_{n'}(\sigma')$ if and only if $n=n'$ and $\sigma \cong\sigma'$. 
\endproclaim 
\subhead 
1.3 
\endsubhead 
We recall the definition of the Artin exponent $\ar(\sigma)$ and the Swan exponent $\sw(\sigma)$, for $\sigma \in \dw F$. 
\par 
Let $\psi$ be a non-trivial smooth character of $F$ and $s$ a complex variable. The Langlands-Deligne local constant $\ve(\sigma,s,\psi)$ takes the form 
$$ 
\ve(\sigma,s,\psi) = \ve(\sigma,0,\psi)\,q^{-(\ar(\sigma)+nc(\psi))s}. 
$$ 
The constant $\ve(\sigma,0,\psi)$ is non-zero. The exponent $\ar(\sigma)$ is a non-negative integer depending only on $\sigma$ and $c(\psi)$ is an integer depending only on $\psi$. The function $\sigma\mapsto \ar(\sigma)$ is additive with respect to direct sums. In simple cases, it is given as follows. 
\proclaim{Fact 1} 
\roster 
\item 
If $\chi$ is an unramified character of $\scr W_F$, then $\ar(\roman{Sp}_r(\chi)) = r{-}1$. 
\item 
Let $(\sigma,\frak n) \in \dw F$. If $\sigma$ is a direct sum of unramified characters of\/ $\scr W_F$, then $\ar(\sigma,\frak n)$ equals the rank of the linear operator $\frak n$. 
\item 
If $\sigma \in \wW F$ is not an unramified character, then $\ar(\roman{Sp}_r(\sigma)) = r\ar(\sigma)$. 
\endroster 
\endproclaim 
These are the key instances of a general formula \cite{23} (4.1.6) (but note that, in the terminology of \cite{23}, all $\sigma\in \dw F$ are $\Phi$-semisimple). 
\par 
We define the {\it Swan exponent\/} $\sw(\sigma)$ of $\sigma\in \dw F$: if $\sigma = (\sigma_\scr W,\frak n)$, then $\sw(\sigma) = \sw(\sigma_\scr W)$. On $\wss F$, the function $\sigma\mapsto \sw(\sigma)$ is additive with respect to direct sums. If $\chi$ is an unramified character of $\scr W_F$, then $\sw(\chi) = \ar(\chi) = 0$. If $\sigma\in \wW F$ is not an unramified character, then $\sw(\sigma) = \ar(\sigma) - \dim\sigma$. 
\par 
As in the introduction, it is helpful to have normalized exponents. For $\sigma\in \dw F$, $\sigma\neq 0$, set  
$$ 
\eta(\sigma) = \ar(\sigma)/\dim\sigma, \qquad \vs(\sigma) = \sw(\sigma)/\dim\sigma.  
$$ 
{\it The use of either of these symbols carries the presumption that $\sigma$ is not zero.} 
\par 
The $\vs$-invariant has a helpful property. For a real number $x\ge 0$, let $\scr W_F^x$ be the corresponding ramification subgroup of $\scr W_F$ (see \cite{20} IV \S3). From \cite{13} Th\'eor\`eme 3.5, we have: 
\proclaim{Fact 2} 
If $\sigma\in \wW F$, then $\vs(\sigma) = \roman{inf}\,\big\{x\ge 0: \scr W_F^x \subset \roman{Ker}\,\sigma\big\}$. 
\endproclaim 
We make repeated use of the following observation. 
\proclaim{Lemma} 
If $\sigma, \tau\in \wW F$, then 
$$ 
\align 
\vs(\sigma\otimes \tau) &\le \roman{max}\,\{\vs(\sigma),\vs(\tau)\}, \quad \text{and}\\ 
\eta(\sigma\otimes\tau) &\le \roman{max}\,\{\eta(\sigma),\eta(\tau)\}. 
\endalign 
$$ 
Equality holds in the first instance if $\vs(\tau) \neq \vs(\sigma)$, in the second if $\eta(\tau)\neq \eta(\sigma)$. 
\endproclaim 
\demo{Proof} 
The first assertion follows from Fact 2. A simple calculation from the definition then gives the second. \qed 
\enddemo 
\remark{Note} 
We have proved the parts of Theorems C and CS relating to irreducible representations. We do not return to those results until section 5.  
\endremark 
\subhead 
1.4 
\endsubhead 
We consider tensor products of indecomposable elements of $\dw F$ that are {\it unramified\/} on restriction to $\scr W_F$. 
\proclaim{Proposition} 
Let $m$, $n$ be positive integers. If $\chi$, $\xi$ are unramified characters of $\scr W_F$ then 
$$ 
\align 
\ar(\roman{Sp}_m(\chi)\otimes \roman{Sp}_n(\xi)) &= mn-\roman{min}\,\{m,n\}, \\ 
\eta(\roman{Sp}_m(\chi)\otimes \roman{Sp}_n(\xi)) &= \roman{max}\,\{\eta(\roman{Sp}_m(\chi)), \eta(\roman{Sp}_n(\xi))\}. 
\endalign 
$$ 
\endproclaim 
\demo{Proof} 
The two assertions are visibly equivalent so we prove the first. There are positive integers $r_i$ and unramified characters $\chi_i$ of $\scr W_F$, $1\le i\le l$, so that 
$$ 
\roman{Sp}_m(\chi)\otimes \roman{Sp}_n(\xi) = \bigoplus_{i=1}^l \roman{Sp}_{r_i}(\chi_i). 
$$ 
In particular, $\sum_{i=1}^l r_i = mn$. Using the definition of $\eta$ and the additivity of the exponent $\ar$, we get 
$$ 
\align 
\eta\big(\tsize\bigoplus_{i=1}^l \roman{Sp}_{r_i}(\chi_i)\big) &= \sum_{i=1}^l r_i\eta\big(\roman{Sp}_{r_i}(\chi_i)\big)/mn \\ &= \sum_{i=1}^l (r_i{-}1)/mn = 1-l/mn. 
\endalign 
$$ 
We therefore need to compute $l$. 
\par 
Write $\roman{Sp}_m(\chi) = (\sigma, \frak m)$, where $\frak m$ is a regular nilpotent endomorphism of $\Bbb C^m$. Likewise write $\roman{Sp}_n(\xi) = (\tau,\frak n)$, so that $\roman{Sp}_m(\chi)\otimes \roman{Sp}_n(\xi) = (\sigma\otimes\tau, \frak l)$, where $\frak l = \frak m\otimes 1 + 1\otimes\frak n$. In this form, the integer $mn{-}l$ is the rank of the nilpotent operator $\frak l$ (1.3 Fact 1(2)). It is therefore enough to recall: 
\proclaim{Lemma} 
Let $\frak m$ (resp\. $\frak n$) be a regular nilpotent endomorphism of the vector space $V = \Bbb C^m$ (resp\. $W = \Bbb C^n$). The operator $\frak l = \frak m\otimes 1_W+1_V\otimes \frak n$ has rank $mn-\roman{min}\{m,n\}$. 
\endproclaim 
The proof of the lemma is a straightforward exercise which completes the proof of the proposition. \qed 
\enddemo 
%%%%%%%%%%%%%%%%%%%%%%%%%%%%%%%%% 
\head\Rm 
2. Irreducible representations 
\endhead 
We prove Theorems A and AS for {\it irreducible\/} representations of $\scr W_F$, taking an indirect approach. We state and prove analogous results for irreducible {\it cuspidal\/} representations of general linear groups $\GL nF$ and then use the Langlands correspondence. 
\subhead 
2.1 
\endsubhead 
We need some definitions. Let $\pi$ be an irreducible {\it cuspidal\/} representation of $\GL nF$, for an integer $n\ge1$. If $\chi$ is a character of $F^\times$, then $\chi\pi$ denotes the representation $g\mapsto \chi(\det g)\pi(g)$, $g\in \GL nF$. 
\par 
We recalled in the introduction the definition of the Artin exponent $\ar (\pi)$ of $\pi$. We also use the notation $\eta(\pi) = \ar(\pi)/n$. 
\par
The {\it Swan exponent\/} $\sw(\pi)$ of $\pi$ is defined by $\sw(\pi) = \ar(\pi){-}n$ {\it except\/} in the case where $n=1$ and $\chi$ is an unramified character of $F^\times = \GL 1F$. In that case, $\sw(\pi) = 0$. In all cases, $\sw(\pi) \ge 0$. We also use the notation $\vs(\pi) = \sw(\pi)/n$. 
\par 
If $\sigma = \upr L\pi \in \wW F$ is the irreducible representation of $\scr W_F$ attached to $\pi$ by the Langlands correspondence, then $\ar(\sigma) = \ar(\pi)$ and $\eta(\sigma) = \eta(\pi)$. The definitions ensure that $\sw(\pi) = \sw(\sigma)$ and $\vs(\pi) = \vs(\sigma)$. 
\par 
We make a similar modification to the Rankin-Selberg exponent $\ar(\pi\times\rho)$ defined in the introduction. 
\proclaim{Definition} 
Let $\rho$ (resp\. $\pi$) be an irreducible cuspidal representation of $\GL mF$ (resp\. $\GL nF$). Let $d$ be the number of unramified characters $\chi$ of $F^\times$ such that $\chi\rho\cong \check\pi$. Set 
$$ 
\align 
\sw(\pi\times\rho) &= \ar(\pi\times\rho) - mn + d, \\ 
\vs(\pi\times\rho) &= \sw(\pi\times\rho)/mn. 
\endalign 
$$ 
\endproclaim  
Note that if, in this definition, we have $d\neq 0$, then $m=n$ and $d$ divides $n$. As a consequence of the definition and corresponding properties of the Artin exponent, we have: 
\proclaim{Fact} 
Let $\pi$ be an irreducible cuspidal representation of $\GL nF$. 
\roster 
\item 
If $\sigma = \upr L\pi \in \wW F$, then $\sw(\pi) = \sw(\sigma)$ and $\vs(\pi) = \vs(\sigma)$. 
\item 
If $\rho$ is an irreducible cuspidal representation of $\GL mF$ and $\upr L\rho = \tau$, then $\sw(\pi\times\rho) = \sw(\sigma\otimes \tau)$ and $\vs(\pi\times\rho) = \vs(\sigma\otimes\tau)$. 
\endroster 
\endproclaim  
\subhead 
2.2 
\endsubhead 
We remark on some upper bounds. 
\proclaim{Proposition} 
For $i=1,2$, let $\pi_i$ be an irreducible cuspidal representation of $\GL{n_i}F$. 
\roster 
\item 
We have $\vs(\pi_1\times\pi_2) \le \roman{max}\,\{\vs(\pi_1), \vs(\pi_2)\}$. If $\vs(\pi_1) \neq \vs(\pi_2)$, then $\vs(\pi_1\times\pi_2) = \roman{max}\,\{\vs(\pi_1), \vs(\pi_2)\}$. 
\item 
We have $\eta(\pi_1\times\pi_2) \le \roman{max}\,\{\eta(\pi_1),\eta(\pi_2)\}$, with equality in the case $\eta(\pi_1) \neq \eta(\pi_2)$. 
\endroster 
\endproclaim 
\demo{Proof} 
This follows from 1.3 Lemma via the Langlands correspondence. \qed 
\enddemo 
\subhead 
2.3 
\endsubhead 
We consider the more substantial problem of lower bounds. 
\proclaim{Proposition} 
For $i=1,2$, let $\pi_i$ be an irreducible cuspidal representation of $\GL{n_i}F$. If $\pi_1$ is $\vs$-minimal then $\vs(\pi_1\times\pi_2) \ge \tfrac12\,\roman{max}\,\{\vs(\pi_1),\vs(\pi_2)\}$. 
\endproclaim 
\demo{Proof} 
By 2.2 Proposition, we have $\vs(\pi_1\times\pi_2) = \roman{max}\,\{\vs(\pi_1),\vs(\pi_2)\}$ provided $\vs(\pi_1) \neq \vs(\pi_2)$. Assume therefore that $\vs(\pi_1) = \vs(\pi_2)$. In that case, $\vs(\pi_1\times\pi_2) \ge \vs(\pi_1\times \check\pi_1)$, by 2.2 Theorem of \cite{2}. We therefore need to prove the following crucial result.  
\proclaim{Lemma} 
Let $\pi$ be an irreducible cuspidal representation of $\GL  nF$, for some $n\ge1$. If $\pi$ is $\vs$-minimal, then $\vs(\pi\times \check\pi) \ge \vs(\pi)/2$. 
\endproclaim 
\demo{Proof} 
If $\vs(\pi) = 0$, then also $\vs(\pi\times\check\pi) = 0$, as an instance of \cite{6} 6.5 Theorem. We therefore assume $\vs(\pi) > 0$. Thus $\pi$ contains a simple character $\theta$, attached to a simple stratum $[\frak a,m,0,\beta]$ in the matrix algebra $\M nF$. The algebra $F[\beta]$ is a field, of degree $d_\beta$, say, over $F$ and ramification index $e_\beta$. Indeed, $\theta$ is ``m-simple'' ({\it cf\.} \cite{4}, especially Corollary 1), so $e_\beta$ equals the $F$-period of the hereditary order $\frak a$ and $\vs(\pi) = m/e_\beta$. The element $\beta$ determines a certain non-negative integer $\frak c(\beta)$, as in \cite{6} 6.4, such that $\vs(\pi\times\check\pi) = \frak c(\beta)/d_\beta^2$. 
\par 
We next choose a simple stratum $[\frak a,m,m{-}1,\alpha]$ equivalent to $[\frak a,m,m{-}1,\beta]$. Let the field extension $F[\alpha]/F$ have degree $d_\alpha$ and ramification index $e_\alpha$. It follows from 3.1 Proposition of \cite{2} that $\frak c(\alpha)/d_\alpha^2 \le \frak c(\beta)/d_\beta^2$. The element $\alpha$ is minimal over $F$ (in the sense of \cite{7} 1.4.14) and, since $\pi$ is $\vs$-minimal, $d_\alpha > 1$. To calculate $\frak c(\alpha)$, we take a simple stratum $[\frak a',m',0,\alpha]$ in the matrix algebra $\End F{F[\alpha]} \cong \M{d_\alpha}F$. The integer $m'$ is $me_\alpha/e_\beta$ and, by 4.1 Proposition of \cite{5}, $\frak c(\alpha) = m'd_\alpha(d_\alpha{-}1))/e_\alpha$. Therefore 
$$ 
\vs(\pi\times\check\pi) \ge \frak c(\alpha)/d_\alpha^2 = (1-d_\alpha^{-1})\,m/e_\beta \ge m/2e_\beta = \tfrac 12\,\vs(\pi), 
$$ 
as required. \qed 
\enddemo 
This completes the proof of the proposition. \qed 
\enddemo 
\example{Example 1} 
Let $\pi$ be an irreducible, cuspidal representation of $\GL 2F$. Suppose that $\pi$ is $\vs$-minimal and $\vs(\pi) > 0$. In the proof of the last lemma, we get $\alpha = \beta$ and $d_\alpha = 2$. This implies $\vs(\pi\times\check\pi) = \frac12 \vs(\pi)$. The constant $\frac12$ in the proposition is therefore best possible as applied to arbitrary representations. 
\endexample 
\example{Example 2} 
One can improve the constant by restricting the class of representations under consideration. For example, if $\ell\ge3$ is a prime number and if $\rho$ is an irreducible, $\vs$-minimal, cuspidal representation of $\GL\ell F$ with $\vs(\rho) > 0$, the same argument gives $\vs(\rho\times\check\rho) = (1{-}\ell^{-1})\vs(\rho)> \frac12\vs(\rho)$. 
\endexample 
We translate in terms of Artin exponents. Let $\pi$ be an irreducible cuspidal representation of $\GL nF$. If $\pi$ is $\eta$-minimal, it is then $\vs$-minimal. (The converse does not hold: the case $n = 1$ and $\ar(\pi) = 1$ provides an example). 
\proclaim{Corollary} 
For $i=1,2$, let $\pi_i$ be an irreducible cuspidal representation of $\GL{n_i}F$. If $\pi_1$ is $\eta$-minimal then $\eta(\pi_1\times\pi_2) \ge \tfrac12\,\roman{max}\,\{\eta(\pi_1),\eta(\pi_2)\}$. 
\endproclaim 
\demo{Proof} 
If either $\pi_i$ is an unramified character of $F^\times$, there is nothing to prove so we assume otherwise. Suppose next that $\pi_2$ is not an unramified twist of $\check\pi_1$. Thus 
$$ 
\align 
\eta(\pi_1\times\pi_2) = \vs(\pi_1\times\pi_2){+}1 &\ge \tfrac12\,\roman{max}\,\{\vs(\pi_1),\vs(\pi_2)\}{+}1 \\ 
&= \tfrac12\,\roman{max}\,\{\eta(\pi_1),\eta(\pi_2)\} + \tfrac12. 
\endalign 
$$ 
Finally suppose $\pi_2$ is an unramified twist of $\check\pi_1$. Thus $\vs(\pi_1\times\pi_2) = \vs(\pi_1\times\check\pi_1)$. The lemma then gives $\vs(\pi_1\times\pi_2) \ge \tfrac12\vs(\pi_1) = \tfrac12\vs(\pi_2)$. In this case, $n_1 = n_2$ and, since $\pi_1$ is not an unramified character of $F^\times$,  we have $n_1>1$. If $d(\pi_1)$ is the number of unramified characters $\chi$ for which $\chi\pi_1 \cong \pi_1$, we have $d(\pi_1)/n_1^2 \le 1/n_1 \le \tfrac12$. So, 
$$ 
\align 
\eta(\pi_1\times\pi_2) &= \vs(\pi_1\times\pi_2){+}1{-}d(\pi_1)/n_1^2 \\ 
&\ge \vs(\pi_1\times\pi_2){+}\tfrac12 \\ 
&\ge \tfrac12\,\roman{max}\,\{\eta(\pi_1),\eta(\pi_2)\}, 
\endalign 
$$ 
as required. \qed 
\enddemo 
\example{Example 3} 
In Example 1, we may choose $\pi$ so that $d(\pi) = 2$. We then get $\eta(\pi\times\check\pi) = \frac12 \eta(\pi)$. The constant $\frac12$ is thus best possible for Artin exponents as well. 
\endexample 
%%%%%%%%%%%%%%%%%%%%%%%%%%%%%%%%%%% 
\head\Rm 
3. First lower bound 
\endhead 
We prove Theorem A, then deal with Theorem AS at the end of the section. 
\subhead 
3.1 
\endsubhead 
We make a simple reduction. 
\proclaim{Proposition} 
Let $\sigma \in \dw F$ be $\eta$-minimal. If the inequality 
$$ 
\eta(\sigma\otimes\tau) \ge \tfrac12\, \roman{max}\{\eta(\sigma), \eta(\tau)\} 
$$ 
holds when $\tau\in \dw F$ is indecomposable, then it holds for all $\tau\in \dw F$. \endproclaim 
\demo{Proof} 
Let $\tau\in \dw F$ and write $\tau = \bigoplus_{j\in J} \tau_j$, where each $\tau_j$ is an indecomposable element of $\dw F$. Put $\alpha_j = \dim\tau_j/\dim\tau$, so that $\sum_{j\in J} \alpha_j = 1$. By hypothesis, $2\eta(\sigma\otimes\tau_j) \ge \eta(\sigma)$, so 
$$ 
2\eta(\sigma\otimes \tau) = \sum_{j\in J} 2\alpha_j\eta(\sigma\otimes\tau_j) \ge \sum_{j\in J} \alpha_j\eta(\sigma) = \eta(\sigma). 
$$ 
The hypothesis also gives $2\eta(\sigma\otimes\tau_j) \ge \eta(\tau_j)$ so 
$$ 
2\eta(\sigma\otimes\tau) \ge \sum_{j\in J} \alpha_j\eta(\tau_j) = \eta (\tau),  
$$ 
as required. \qed 
\enddemo 
We therefore have to prove: 
\proclaim{Theorem} 
If $\sigma \in \dw F$ is $\eta$-minimal and $\tau\in \dw F$ is indecomposable, then 
$$ 
\eta(\sigma\otimes\tau) \ge \tfrac12\, \roman{max}\{\eta(\sigma), \eta(\tau)\}. 
$$ 
\endproclaim 
This will take us to the end of 3.4. 
\subhead 
3.2 
\endsubhead 
Let $\sigma\in \dw F$. Say that $\sigma$ is {\it $\eta$-homogeneous\/} if there exists $a\in \Bbb R$ such that $\eta(\tau) = a$, for every irreducible factor $\tau$ of $\sigma$ on $\scr W_F$. When this holds, we write $\ell_0(\sigma) = a$. 
\example{Example} 
If $\tau$ is irreducible and $\sigma = \roman{Sp}_r(\tau)$, then $\sigma$ is $\eta$-homogene\-ous with $\ell_0(\sigma) = \eta(\tau)$. Moreover, $\ell_0(\sigma) = \eta(\sigma)$ if $\ell_0(\sigma) \neq 0$. If $\ell_0(\sigma) = 0$, then $\tau$ is an unramified character and $\eta(\sigma) = 1{-}r^{-1}$. 
\endexample 
\proclaim{Proposition} 
Let $\sigma\in \dw F$ be $\eta$-homogeneous and write $\sigma = \bigoplus_{j\in J} \roman{Sp}_{r_j}(\tau_j)$, where the $\tau_j$ are irreducible and $r_j\ge 1$. 
\roster 
\item 
The representation $\sigma$ is $\eta$-minimal if and only if each $\tau_j$ is $\eta$-minimal. 
\item 
If $\sigma$ is $\eta$-minimal, then 
$$ 
\eta(\chi\otimes\sigma) = \roman{max}\,\{\eta(\sigma), \eta(\chi)\},  
$$ 
for all characters $\chi$ of\/ $\scr W_F$. 
\endroster 
\endproclaim 
\demo{Proof} 
Set $a = \ell_0(\sigma)$. If $a=0$, then $\sigma$ is $\eta$-minimal, each $\roman{Sp}_{r_j}(\tau_j)$ is $\eta$-minimal and (2) is immediate. We assume henceforth that $a > 0$. 
\par 
Take first the case where $\sigma$ is {\it indecomposable,\/} say $\sigma = \roman{Sp}_r(\tau)$. If $\dim\tau = 1$, it is clear that $\sigma$ is $\eta$-minimal if and only if $\eta(\tau)  = 0$, contrary to hypothesis. Therefore $\dim\tau > 1$ and $a = \eta(\sigma) = \eta(\tau)$. Let $\chi$ be a character of $\scr W_F$, and set $\eta(\chi) = c$. If $c > a$, then 
$$ 
\eta(\chi\otimes\sigma) = \eta(\chi\otimes\tau) = c > a = \eta(\sigma), 
$$ 
(1.3 Lemma). If, however, $c<a$, we get $\eta(\chi\otimes\sigma) = a = \eta(\sigma) > c$. Suppose finally that  $c=a$. Since $\dim\tau > 1$, we get $\eta(\chi\otimes\sigma) = \eta(\chi\otimes\tau)$. If $\tau$ is not $\eta$-minimal, we may choose $\chi$ so that $\eta(\chi\otimes\tau) < \eta(\tau)$ and so $\sigma$ is not $\eta$-minimal. If $\tau$ is $\eta$-minimal, $\eta(\chi\otimes\tau) = \eta(\tau)$ and we are done with the case of $\sigma$ indecomposable. 
\par 
For the general case, we set $\sigma = \bigoplus_{j\in J} \sigma_j$, where $\sigma_j$ is indecomposable. Put $\alpha_j = \dim\sigma_j/\dim\sigma$, so that 
$$ 
\sum_j\alpha_j = 1 \quad \text{and} \quad \eta(\chi\otimes\sigma) = \sum_j\alpha_j\eta(\chi\otimes\sigma_j). 
$$ 
Suppose $\sigma_j$ is not $\eta$-minimal, for some $j\in J$. We have just shown that there exists $\chi$ with $\eta(\chi) = a$ and $\eta(\chi\otimes\sigma_j) < \eta(\sigma_j)$. On the other hand, $\eta(\chi\otimes\sigma_k) \le \eta(\sigma_k)$ for $k\neq j$ ({\it cf\.} 1.3 Lemma), so 
$$ 
\eta(\chi\otimes\sigma) = \sum_{i\in J} \alpha_i\eta(\chi\otimes\sigma_i) < \sum_{i\in J} \alpha_i\eta(\sigma_j) = \eta(\sigma), 
$$ 
whence $\sigma$ is not $\eta$-minimal. 
\par 
Assume, therefore, that every $\sigma_j$  is $\eta$-minimal. Let $c = \eta(\chi)$. If $c\ge a$, the discussion of the indecomposable case gives $\eta(\chi\otimes\sigma_j) = c$, $j\in J$, so $\eta(\chi\otimes\sigma) = c \ge \eta(\sigma)$. If, however, $c<a$, we get $\eta(\chi\otimes\sigma) = a = \eta(\sigma)$. Thus $\sigma$ is $\eta$-minimal and we have also proven (2). \qed 
\enddemo 
\subhead 
3.3 
\endsubhead 
We use 3.2 Proposition to prove a special case of 3.1 Theorem. 
\proclaim{Proposition} 
Let $\sigma\in \dw F$ be $\eta$-minimal and $\eta$-homogeneous. If $\tau \in \dw F$ is irreducible, then 
$$ 
\eta(\sigma\otimes\tau) \ge \tfrac12\,\roman{max}\,\{\eta(\sigma),\eta(\tau)\}. 
$$ 
\endproclaim 
\demo{Proof} 
Combining 3.2 Proposition with the decomposition technique used in 3.1, we reduce to the case where $\sigma$ is indecomposable and $\eta$-minimal. Let $a = \ell_0(\sigma)$.
\par
Consider first the case where $a = 0$, that is, $\sigma = \roman{Sp}_r(\chi)$ with $\chi$ an unramified character of $\scr W_F$ and $r\ge 1$. Thus $\eta(\sigma) = (r{-}1)/r$. On the other hand, $\sigma\otimes\tau = \roman{Sp}_r(\chi\otimes\tau)$ and so 
$$ 
\eta(\sigma\otimes\tau) = \left\{\,\alignedat2 &\eta(\tau) \quad &\text{if $\eta(\tau)\neq 0$,} \\ 
&\eta(\sigma) \quad &\text{if $\eta(\tau) = 0$.} \endalignedat \right. 
$$ 
In the first case, we have $\eta(\tau) \ge 1 > \eta(\sigma)$ while, in the second, $\eta(\tau) \le \eta(\sigma)$. The result therefore holds when $a=0$. 
\par 
From now on, we assume $a>0$. We write $\sigma$ as $\roman{Sp}_r(\rho)$, where $\rho$ is irreducible and $\eta$-minimal. We have $\eta(\rho) = \eta(\sigma) = a$ and $\dim\rho > 1$. Assume initially that there is no unramified character $\chi$ of $\scr W_F$ such that $\chi\otimes \rho \cong \check\tau$. This means that no irreducible component of $\rho\otimes \tau$ is unramified, so $\ar(\sigma\otimes \tau) = r\ar(\rho\otimes\tau)$ and 
$$ 
2\eta(\sigma\otimes \tau) = 2\eta(\rho\otimes\tau) \ge \roman{max}\,\{\eta(\rho),\eta(\tau)\}, 
$$ 
by 2.3 Corollary. Since $\eta(\sigma) = \eta(\rho)$, we are done in this case. 
\par 
For the remaining case, we may assume $\check\tau \cong \rho$: in particular, $\eta(\tau) = a$. Let $d$ be the number of unramified characters $\chi$ for which $\chi\otimes\rho\cong \rho$. Thus $d$ divides $m = \dim\rho > 1$. To estimate $\ar(\sigma\otimes\tau) = \ar(\roman{Sp}_r(1)\otimes \rho\otimes\check\rho)$, we write 
$$ 
\rho\otimes\check\rho = \rho'\oplus \chi_1\oplus \dots \oplus \chi_d, 
$$ 
where the $\chi_i$ are unramified characters of $\scr W_F$ and every irreducible component of $\rho'$ has strictly positive exponent. Thus 
$$ 
\eta(\rho\otimes\check\rho) = (1-d/m^2)\eta(\rho'). 
$$ 
Also, $\eta(\rho\otimes\check\rho) \ge \tfrac12 \eta(\rho)$ by 2.3 Corollary. Taking this into account, we have 
$$ 
\align 
\eta(\sigma\otimes\tau) &= \eta\big((\roman{Sp}_r(1)\otimes \rho') \oplus \sum_{i=1}^d\roman{Sp}_r(\chi_i)\big) \\ 
&= (1{-}d/m^2)\,\eta(\roman{Sp}_r(1)\otimes \rho') + d(r{-}1)/rm^2 \\ 
&= (1{-}d/m^2)\,\eta(\rho') + d(r{-}1)/rm^2 \\ 
&\ge (1{-}d/m^2)\,\eta(\rho') \\ 
&\ge \tfrac12 \eta(\rho) = \tfrac12\eta(\tau).
\endalign 
$$ 
Since, in this case, $\eta(\rho) = \eta(\sigma)$ the proof is complete. \qed 
\enddemo 
We may now deal with 3.1 Theorem in the case where $\sigma$ is $\eta$-homogeneous. 
\proclaim{Corollary} 
Let $\sigma\in \dw F$ be $\eta$-minimal and $\eta$-homogeneous. If $\tau \in \dw F$ is indecomposable, then $\eta(\sigma\otimes\tau) \ge \frac12\,\roman{max}\,\{\eta(\sigma),\eta(\tau)\}$. 
\endproclaim 
\demo{Proof} 
As in the proof of the proposition, it is enough to treat the case where $\sigma$ is indecomposable and $\eta$-minimal. Thus $\sigma = \roman{Sp}_r(\sigma')$ and $\tau = \roman{Sp}_s(\tau')$, for integers $r,s \ge 1$ and irreducible representations $\sigma'$, $\tau'$. We have 
$$ 
\sigma\otimes\tau = \big(\roman{Sp}_r(1)\otimes \roman{Sp}_s(1) \otimes \sigma'\big) \otimes\tau'. 
$$ 
The representation $\sigma'' = \roman{Sp}_r(1)\otimes \roman{Sp}_s(1) \otimes \sigma'$ is $\eta$-minimal and $\eta$-homogeneous with $\ell_0(\sigma'') = \ell_0(\sigma) = \eta(\sigma')$, so the proposition gives 
$$ 
\eta(\sigma\otimes \tau) = \eta(\sigma''\otimes\tau') \ge \tfrac12\,\roman{max}\,\{\eta(\sigma''),\eta(\tau')\}. 
$$ 
It is therefore enough to show that 
$$ 
\roman{max}\,\{\eta(\sigma''),\eta(\tau')\} \ge \roman{max}\,\{\eta(\sigma),\eta(\tau)\}. 
\tag $*$ 
$$ 
\indent  
To do this, we write the tensor product $\roman{Sp}_r(1)\otimes \roman{Sp}_s(1)$ as a sum of indecomposable representations: there are unramified characters $\chi_i$ and positive integers $r_i$, $1\le i \le l$, such that 
$$ 
\roman{Sp}_r(1)\otimes \roman{Sp}_s(1) = \bigoplus_{i=1}^l \roman{Sp}_{r_i}(\chi_i). 
$$ 
We have $\sum_{i=1}^l r_i = rs$ and, by 1.4 Proposition, $l = \roman{min}\,\{r,s\}$. Accordingly, 
$$ 
\sigma'' = \bigoplus_{i=1}^l \roman{Sp}_{r_i}(\chi_i\otimes\sigma'). 
$$ 
If $\sigma'$ is not unramified, then $\eta(\sigma'') = \eta(\sigma') = \eta(\sigma)$. Likewise, if $\tau$ is not unramified then $\eta(\tau') = \eta(\tau)$. So, if neither $\sigma'$ nor $\tau'$ is unramified, we get $\roman{max}\,\{\eta(\sigma''),\eta(\tau')\} = \roman{max}\,\{\eta(\sigma),\eta(\tau)\}$, proving $(*)$ in this case.  
\par 
Suppose next that $\sigma'$ is unramified. By 1.4 Proposition, we have 
$$ 
\eta(\sigma'') = \sum_i (r_i{-}1)/rs = 1-\tfrac l{rs}, 
$$ 
while $\eta(\sigma) = (r{-}1)/r \le \eta(\sigma'')$. If $\tau'$ is not unramified, then $\eta(\tau') = \eta(\tau)$ and we are done. If $\tau'$ is unramified, then $\eta(\tau') = 0$ and $\eta(\tau) = (s{-}1)/s \le \eta(\sigma'')$, since $l\le s$.  This proves $(*)$ in all cases, and the proof is complete. \qed 
\enddemo 
\subhead 
3.4 
\endsubhead 
We enter the final stage of the proof of 3.1 Theorem. We proceed in two steps.  
\proclaim{Proposition 1} 
If $\sigma\in \dw F$ is $\eta$-minimal  then $\eta(\sigma\otimes\tau) \ge \tfrac12\,\eta(\sigma)$, for all indecomposable $\tau\in \dw F$. 
\endproclaim 
\demo{Proof} 
It is enough to prove that $\eta(\xi_i\otimes\tau) \ge \tfrac12\,\eta(\xi_i)$, for each indecomposable component $\xi_i$ of $\sigma$. There is a character $\chi$ such that $\chi\otimes\tau$ is $\eta$-minimal. Since $\tau$ is indecomposable, $\chi\otimes\tau$ is $\eta$-homogeneous so 3.3 Corollary gives 
$$ 
\eta(\xi_i \otimes \tau) = \eta\big((\chi^{-1}\otimes\xi_i) \otimes (\chi\otimes\tau)\big) \ge \tfrac12\,\eta(\chi^{-1}\otimes\xi_i). 
$$ 
Consequently, 
$$ 
\eta(\sigma\otimes\tau) = \eta(\chi^{-1}\otimes\sigma \otimes\chi\otimes\tau) \ge \tfrac12 \eta(\chi^{-1}\otimes\sigma). 
$$ 
As $\sigma$ is $\eta$-minimal, so $\eta(\chi^{-1}\otimes\sigma) \ge \eta(\sigma)$ and the result follows. \qed 
\enddemo 
It now remains only to prove: 
\proclaim{Proposition 2} 
If $\sigma\in \dw F$ is $\eta$-minimal and $\tau\in \dw F$ is indecomposable, then $\eta(\sigma\otimes\tau) \ge \tfrac12\,\eta(\tau)$. 
\endproclaim 
\demo{Proof} 
The representation $\tau$ is $\eta$-homogeneous. If $\tau$ is $\eta$-minimal the result follows from 3.3 Corollary and 3.1 Proposition. We therefore assume the contrary. 
\par 
To proceed further, we need to extend 3.2 Proposition. Write 
$$ 
\sigma = \bigoplus_{i\in I} \roman{Sp}_{r_i}(\xi_i), 
$$ 
for irreducible representations $\xi_i$ and integers $r_i\ge 1$. Let $c = \roman{max}\,\eta(\xi_i)$. Define $\sigma_\roman{max}$ as the sum of all factors $\roman{Sp}_{r_i}(\xi_i)$ for which $\eta(\xi_i) = c$, and $\sigma'$ as the sum of the others. 
\proclaim{Lemma} 
The representation $\sigma'$ is either zero or $\eta$-minimal. 
\endproclaim 
\demo{Proof}
Assume $\sigma'\neq 0$. Let $d_\roman{max} = \dim\sigma_\roman{max}$ and $d' = \dim\sigma'$. Set $d = d_\roman{max}{+}d' = \dim\sigma$. Let $\phi$ be a character of $\scr W_F$ and write $s = \eta(\phi)$. We compare the expressions 
$$ 
\align 
d\eta(\sigma) &= d_\roman{max}\eta(\sigma_\roman{max}) + d'\eta(\sigma'), \\ 
d\eta(\phi\otimes\sigma) &= d_\roman{max}\eta(\phi\otimes \sigma_\roman{max}) + d'\eta(\phi\otimes\sigma'). 
\endalign 
$$ 
If $s<c$, then $\eta(\phi\otimes\sigma_\roman{max}) = \eta(\sigma_\roman{max}) = c$. Since $\eta(\phi\otimes\sigma) \ge \eta(\sigma)$, we have $\eta(\phi\otimes\sigma') \ge \eta(\sigma')$. On the other hand, if $s\ge c$, then $\eta(\phi\otimes\sigma') = s \ge c > \eta(\sigma')$. This shows that $\eta(\phi\otimes\sigma') \ge \eta(\sigma')$ for all $\phi$. \qed 
\enddemo 
If the representation $\sigma'$ is zero, then $\sigma = \sigma_\roman{max}$ is $\eta$-homogeneous and $\eta$-minimal. The proposition in this case is given by 3.3 Corollary. We assume therefore that $\sigma' \neq 0$. Certainly $\sigma_\roman{max} \neq 0$ so, using induction on the Jordan-H\"older length of $\sigma$, we may assume 
$$ 
\eta(\sigma'\otimes\tau) \ge \tfrac12\,\eta(\tau). 
$$ 
Let $s = \eta(\tau)$ and let $\chi$ be a character such that $\chi\otimes\tau$ is $\eta$-minimal. Since we assume $\tau$ is not $\eta$-minimal, we have $\eta(\chi) = s > 0$ (1.3 Lemma). 
\par 
Examining cases, suppose first that $s>c$. Here, $\eta(\sigma\otimes\tau) = s >s/2$, implying the result in this situation. If, on the other hand, $s<c$, we have $\eta(\sigma_\roman{max}\otimes\tau) = c$ while, by inductive hypothesis, $\eta(\sigma'\otimes\tau) \ge s/2$. Writing $\alpha = d_\roman{max}/d$, $\beta = d'/d$, we have 
$$ 
\align 
\eta(\sigma\otimes\tau) &= \alpha\eta(\sigma_\roman{max}\otimes \tau) + \beta \eta(\sigma'\otimes\tau) \\ 
&\ge \alpha c+\beta s/2 \ge s/2. 
\endalign 
$$ 
It remains to treat the case $s=c$. Here, $\eta(\sigma'\otimes\tau) = s = c$, so 
$$ 
\eta(\sigma\otimes\tau) = \alpha\eta(\sigma_\roman{max}\otimes\tau) + \beta s. 
$$ 
However, $\eta(\sigma_\roman{max}\otimes\tau) = \eta\big((\chi^{-1}\otimes\sigma_\roman{max})\otimes(\chi\otimes\tau)\big)$ while, by the very first case of this proof, we have 
$$ 
\eta\big((\chi^{-1}\otimes\sigma_\roman{max})\otimes(\chi\otimes\tau)\big) \ge \tfrac12\,\eta(\chi^{-1}\otimes\sigma_\roman{max}). 
$$ 
Since $\sigma$ is $\eta$-minimal, 
$$ 
\align 
\eta(\chi^{-1}\otimes\sigma) &= \alpha\eta(\chi^{-1}\otimes\sigma_\roman{max}) + \beta \eta(\chi^{-1}\otimes\sigma') \\ 
&= \alpha\eta(\chi^{-1}\otimes\sigma_\roman{max}) + \beta s \\ 
&\ge \eta(\sigma) = \alpha s+\beta\eta(\sigma'). 
\intertext{That is,} 
\alpha\eta(\chi^{-1} \otimes \sigma_\roman{max}) &\ge (\alpha{-}\beta)s + \beta\eta(\sigma') 
\intertext{and, overall, }
\eta(\sigma\otimes\tau) &= \alpha\eta\big((\chi^{-1}\otimes\sigma_\roman{max})\otimes(\chi\otimes \tau)\big) + \beta s \\ 
&\ge \tfrac12\alpha\eta(\chi^{-1}\otimes\sigma_\roman{max}) + \beta s \\ 
&\ge \tfrac12 (\alpha{-}\beta)s +\tfrac12\beta\eta(\sigma')+\beta s \\ 
&\ge \tfrac12(\alpha{+}\beta)s \\
&= s = \tfrac12 \eta(\tau). 
\endalign 
$$ 
That is, $ \eta(\sigma\otimes\tau) \ge \frac12 \eta(\tau)$ as required. \qed 
\enddemo 
This completes the proofs of 3.1 Theorem and Theorem A. \qed 
\subhead 
3.5 
\endsubhead 
We digress to highlight a special case. Say that $\sigma\in \dw F$ is {\it unramified\/} if its restriction to $\scr W_F$ is a sum of unramified characters. Any such $\sigma$ is both $\eta$-minimal and $\eta$-homogeneous. 
\proclaim{Example 1} 
If $\sigma\in \dw F$ is unramified, then $\eta(\sigma\otimes\tau) \ge \roman{max}\{\eta(\sigma),\eta(\tau)\}$, 
for all $\tau\in \dw F$. 
\endproclaim 
To justify this, one applies the argument of 3.1 twice to reduce to the case where both $\sigma$ and $\tau$ are indecomposable. The proof of 3.3 Proposition gives the result when $\tau$ is irreducible. In the proof of 3.3 Corollary, we still get $\eta(\sigma\otimes\tau) = \eta(\sigma''\otimes\tau')$, so $\eta(\sigma\otimes\tau) \ge \roman{max}\,\{\eta(\sigma''),\eta(\tau')\}$. In the same proof, we have shown that $\roman{max}\,\{\eta(\sigma''),\eta(\tau')\} \ge \roman{max}\,\{\eta(\sigma),\eta(\tau)\}$, whence the assertion. 
\par 
Starting again from the first case of the proof of 3.3 Proposition, one may equally conclude: 
\proclaim{Example 2} 
Let $\sigma\in \dw F$ be unramified. If $\tau\in \dw F$ has no unramified direct factor, then $\eta(\sigma\otimes \tau) = \eta(\tau)$. 
\endproclaim 
For a tensor product of unramified representations, one may derive an explicit formula from 1.4 Lemma. 
\subhead 
3.6 
\endsubhead 
We prove Theorem AS: if $\sigma,\tau\in \wss F$ and if $\sigma$ is $\vs$-minimal, then $\vs(\sigma\otimes \tau) \ge \frac12\,\roman{max}\{\vs(\sigma),\vs(\tau)\}$. 
\par 
If $\sigma$ and $\tau$ are irreducible, this follows from 2.3 Proposition. An argument identical to 3.1 Proposition shows it is enough to prove the theorem under the additional hypothesis that $\tau$ is irreducible. 
\par 
Say that $\sigma\in \wss F$ is $\vs$-homogeneous if, for some $a$, we have $\vs(\sigma') = a$ for all irreducible components $\sigma'$ of $\sigma$. With this definition, the analogue of 3.2 Proposition holds with the same proof. In light of the case already done, where $\sigma$ and $\tau$ are irreducible, the analogue of 3.3 Proposition is immediate here and the Corollary is redundant. The propositions of 3.4 hold, with identical proofs, and the theorem is proved. \qed 
%%%%%%%%%%%%%%%%%%%%%%%%%%%%%%%%%% 
\head\Rm 
4. Symmetric lower bound 
\endhead 
We prove Theorem B and deal with Theorem BS at the end of the section. We first accumulate some preliminary results concerning irreducible or indecomposable representations. 
\subhead 
4.1 
\endsubhead 
We start with what amounts to a special case of the theorem. 
\proclaim{Proposition} 
If $\sigma$, $\tau$ are irreducible representations of\/ $\scr W_F$, then 
$$ 
\align 
\vs(\sigma\otimes\check\tau) &\ge \roman{max}\,\{\vs(\sigma\otimes\check\sigma), \vs(\tau\otimes\check\tau)\}, \quad \text{and}\\
\eta(\sigma\otimes\check\tau) &\ge \roman{max}\,\{\eta(\sigma\otimes\check\sigma), \eta(\tau\otimes\check\tau)\}. 
\endalign 
$$ 
\endproclaim 
\demo{Proof} 
For the first assertion, we follow \cite{12} but use the notation and layout of \cite{5} 2.5, 3.1. The set $\wW F$ carries a canonical pairing $\Delta$ with non-negative real values \cite{5} (2.5.3). It has the property $\Delta(\sigma,\sigma) \le \Delta(\sigma,\tau)$, for all $\tau \in \wW F$. As in \cite{5} 3.1, there is a continuous, strictly increasing function $\vS_\sigma$ such that $\vS_\sigma(\Delta(\sigma,\tau)) = \vs(\sigma\otimes\check\tau)$. Therefore $\vs(\sigma\otimes\check\sigma) \le \vs(\sigma\otimes\check\tau)$, as desired. 
\par 
In the second assertion, suppose first that $\sigma\not\cong\chi\otimes \tau$, for any unramified character $\chi$ of $\scr W_F$. It follows that $\eta(\sigma\otimes\check\tau) = \vs(\sigma\otimes\check\tau){+}1$. The first assertion then gives $\eta(\sigma\otimes\check\tau) \ge \roman{max}\,\{\vs(\sigma\otimes\check\sigma){+}1, \vs(\tau\otimes\check\tau){+}1\}$. However, 
$$ 
\eta(\sigma\otimes\check\sigma) = \vs(\sigma\otimes\check\sigma)+1-d_\sigma/m^2, 
$$ 
where $d_\sigma$ is the number of unramified characters $\chi$ such that $\chi\otimes\sigma\cong \sigma$ and $m = \dim\sigma$. Likewise for $\tau$, and the result follows. 
\par
If, on the other hand, there is an unramified character $\phi$ such that $\tau\cong \phi\otimes\sigma$, we get $\eta(\sigma\otimes\check\tau) = \eta(\sigma\otimes\check\sigma) = \eta(\tau\otimes\check\tau)$, and there is nothing to do. \qed 
\enddemo 
\subhead
4.2 
\endsubhead 
The exponent has a striking ultrametric property. 
\proclaim{Proposition} 
If $\sigma,\tau, \rho\in \wW F$, then 
$$ 
\align 
\vs(\sigma\otimes\check\tau) &\le \roman{max}\,\{\vs(\sigma\otimes\check\rho), \vs(\rho\otimes\check\tau)\}, \\ 
\eta(\sigma\otimes\check\tau) &\le \roman{max}\,\{\eta(\sigma\otimes\check\rho), \eta(\rho\otimes\check\tau)\}
\endalign 
$$ 
\endproclaim 
\demo{Proof} 
The first inequality is 3.1 Corollary of \cite{5}. To deduce the second, let $d_{\sigma\tau}$ be the number of unramified characters $\chi$ for which $\tau\cong \chi\otimes\sigma$, and similarly for the other pairs. Let $m = \dim\sigma$, $n = \dim\tau$ and $l = \dim\rho$. Thus 
$$ 
\eta(\sigma\otimes\check\tau) = \vs(\sigma\otimes\check\tau){+}1-d_{\sigma\tau}/mn, 
$$ 
and similarly for the others. The first part of the proposition yields 
$$ 
\eta(\sigma\otimes\check\tau) \le \roman{max}\,\{\eta(\sigma\otimes\check\rho)+d_{\sigma\rho}/ml, \eta(\rho\otimes\check\tau) + d_{\rho\tau}/nl\} - d_{\sigma\tau}/mn. 
$$ 
This gives the result if $d_{\sigma\tau} = d_{\sigma\rho} = d_{\rho\tau} = 0$. If $d_{\sigma\tau} \neq 0$, then $m=n$ and $d_{\sigma\rho} = d_{\rho\tau}$. Also, $\eta(\sigma\otimes\check\tau) = \eta(\sigma\otimes\check\sigma)$ and $\eta(\rho\otimes\check\tau) = \eta(\rho\otimes\check\sigma)$. The desired inequality thus reduces to $\eta(\sigma\otimes\check\sigma) \le \eta(\sigma\otimes\check\rho)$, which follows from 4.1 Proposition. Similarly, if $d_{\sigma\rho}\neq 0$, we have to check that $\eta(\sigma\otimes\check\tau) \le \roman{max}\,\{\eta(\sigma\otimes\check\sigma),\eta(\sigma\otimes\check\tau)\}$, and this is immediate. \qed 
\enddemo 
\subhead 
4.3 
\endsubhead 
We generalize the propositions of 4.1, 4.2 to indecomposable representations. To do this, we need some explicit formulas. 
\par 
Let $\sigma$, $\tau$ be irreducible representations of $\scr W_F$ of dimension $m$, $n$ respectively. Let $d_\sigma$ be the number of unramified characters $\chi$ such that $\sigma\cong \chi\otimes\sigma$. Define $d_\tau$ similarly, and let $d_{\sigma\tau}$ be the number of unramified characters $\chi$ such that $\sigma\cong \chi\otimes \tau$. Let $r\ge s\ge 1$ integers, and set $\Sigma = \roman{Sp}_r(\sigma)$, $\Tau = \roman{Sp}_s(\tau)$. 
\proclaim{Lemma} 
With the notation above, we have 
$$ 
\align 
\eta(\Sigma\otimes\check\Tau) &= \eta(\sigma\otimes\check\tau) + d_{\sigma\tau}(1{-}r^{-1})/mn, \\
\eta(\Sigma\otimes\check\Sigma) &= \eta(\sigma\otimes\check\sigma) + d_{\sigma}(1{-}r^{-1})/m^2, \\ 
\eta(\Tau\otimes\check\Tau) &= \eta(\tau\otimes\check\tau) + d_{\tau}(1{-}s^{-1})/n^2. 
\endalign 
$$ 
\endproclaim 
\demo{Proof} 
The second and third relations are instances of the first, so we need only prove that one. 
\par 
We write $\sigma\otimes\check\tau = \rho \oplus \chi_1\oplus\chi_2\oplus \dots\oplus \chi_d$, where every component of $\rho$ is not unramified and $\chi_j$ is an unramified character, $1\le j\le d = d_{\sigma\tau}$. Thus 
$$ 
\eta(\sigma\otimes\check\tau) = (mn{-}d)\eta(\rho)/mn . 
$$ 
We also have 
$$ 
\Sigma\otimes \check\Tau = \roman{Sp}_r(1)\otimes \roman{Sp}_s(1) \otimes \sigma \otimes \check\tau. 
$$ 
Set $\Rho = \roman{Sp}_r(1)\otimes \roman{Sp}_s(1) \otimes\rho$, so that $\eta(\Rho) = \eta(\rho)$. Expanding, we get 
$$ 
\eta(\Sigma\otimes\check\Tau) = (mn{-}d)\eta(\rho)/mn + d\eta(\roman{Sp}_r(1)\otimes \roman{Sp}_s(1) \otimes \vS_j\chi_j)/mn. 
$$ 
We use 1.4 Proposition (and the hypothesis $r\ge s$) to get 
$$ 
\roman{Sp}_r(1)\otimes \roman{Sp}_s(1) = \bigoplus_{k=1}^s \roman{Sp}_{r_k}(\xi_k), 
$$ 
for positive integers $r_k$, with sum $rs$, and unramified characters $\xi_k$. In particular, 
$$ 
\eta(\roman{Sp}_r(1)\otimes \roman{Sp}_s(1) \otimes \phi) = 1{-}r^{-1}, 
$$ 
for any unramified character $\phi$. Therefore 
$$ 
\eta(\roman{Sp}_r(1)\otimes \roman{Sp}_s(1) \otimes \vS_j\chi_j) = d^{-1}\sum_j \eta(\roman{Sp}_r(1)\otimes \roman{Sp}_s(1) \otimes\chi_j) = 1{-}r^{-1} 
$$ 
and, altogether, 
$$ 
\align 
\eta(\Sigma\otimes\check\Tau) &= (mn{-}d)\eta(\rho)/mn + d(1{-}r^{-1})/mn \\ 
&= \eta(\sigma\otimes\check\tau) + d(1{-}r^{-1})/mn, 
\endalign 
$$ 
as required. \qed 
\enddemo 
\proclaim{Proposition 1} 
If $\Sigma,\Tau\in \dw F$ are indecomposable, then 
$$ 
\eta(\Sigma\otimes\check\Tau) \ge \roman{max}\,\{\eta(\Sigma\otimes\check\Sigma), \eta(\Tau\otimes\check\Tau)\}. 
$$ 
\endproclaim 
\demo{Proof} 
We write $\Sigma = \roman{Sp}_r(\sigma)$ and $\Tau = \roman{Sp}_s(\tau)$, for $\sigma,\tau \in \wW F$. Using the notation of the lemma, suppose first that $d_{\sigma\tau} = 0$. Using 4.1 Proposition and the formulas from the lemma, we get 
$$ 
\align 
\eta(\Sigma\otimes\check\Tau) = \eta(\sigma\otimes\check\tau) &= \vs(\sigma\otimes\check\tau){+}1 \\ 
&\ge \vs(\sigma\otimes\check\sigma){+}1 = \eta(\sigma\otimes\check\sigma) + d_\sigma/m^2 \ge \eta(\Sigma\otimes\check\Sigma), 
\endalign 
$$  
and likewise $\eta(\Sigma\otimes\check\Tau) \ge \eta(\Tau\otimes\check\Tau)$. 
\par 
Suppose therefore that $d_{\sigma\tau} \neq 0$. There is then an unramified character $\chi$ for which $\tau\cong \chi\otimes \sigma$. In particular, $m=n$ and $d_{\sigma\tau} = d_\sigma = d_\tau = d$, say. From the formulas above, we get $\eta(\Sigma\otimes\check\Tau) = \eta(\Sigma\otimes\Sigma) \ge \eta(\Tau\otimes\check\Tau)$, which is enough. \qed 
\enddemo 
\proclaim{Proposition 2} 
If $\Sigma,\Tau, \Rho \in \dw F$ are indecomposable then 
$$ 
\eta(\Sigma\otimes\check\Tau) \le \roman{max}\{\eta(\Sigma\otimes \check\Rho), \eta(\Rho\otimes \check\Tau)\}. 
$$ 
\endproclaim 
\demo{Proof} 
There are representations $\sigma,\tau,\rho\in \wW F$ and integers $r$, $s$, $t$ such that $\Sigma = \roman{Sp}_r(\sigma)$, $\Tau = \roman{Sp}_s(\tau)$ and $\Rho = \roman{Sp}_t(\rho)$. Let $\dim\sigma = m$, $\dim\tau = n$ and $\dim\rho = l$. Define integers $d_{\sigma\tau}$, $d_\sigma$ etc., as before. 
\par 
Take first the case where $\sigma$ is not an unramified twist of $\tau$. That is, $d_{\sigma\tau} = 0$ and $\eta(\Sigma\otimes\check\Tau) = \eta(\sigma\otimes\check\tau)$. If, for example, $d_{\sigma\rho}\neq 0$ then $l=m$ and 
$$ 
\eta(\Sigma\otimes\check\Rho) = \eta(\sigma\otimes\check\rho) + d_{\sigma\rho}(1{-}q^{-1})/m^2 \ge \eta(\sigma\otimes \check\rho), 
$$ 
where $q = \roman{max}\{r,t\}$. If, on the other hand, $d_{\sigma\rho} = 0$, we get the conclusion $\eta(\Sigma\otimes\check\Rho) = \eta(\sigma\otimes\check\rho)$. Similarly for the pair $(\rho,\tau)$, so the desired inequality now follows from 4.2 Proposition. 
\par 
We therefore assume $d_{\sigma\tau} \neq 0$. Our assumption $r\ge s$ implies $\eta(\Sigma\otimes\check\Tau) = \eta(\Sigma\otimes\check\Sigma)$, while $\eta(\Sigma\otimes\check\Sigma) \le \eta(\Sigma\otimes\check\Rho)$, by Proposition 1 above. \qed 
\enddemo 
\subhead 
4.4 
\endsubhead 
We now prove the main statement of Theorem B, that is: 
\proclaim{Theorem} 
If $\sigma_1,\sigma_2\in \dw F$, then 
$$ 
\eta(\sigma_1\otimes\sigma_2) \ge \tfrac12\,\big(\eta(\sigma_1\otimes\check\sigma_1)+ \eta(\sigma_2\otimes \check\sigma_2)\big). 
$$ 
\endproclaim 
\demo{Proof} 
We proceed by induction on $r_1r_2$, where $r_i$ is the number of isomorphism classes of indecomposable direct factors of $\sigma_i$. The case $r_1r_2 = 1$ follows from 4.3 Proposition 1, so we assume $r_1r_2 \ge2$. 
\par
If $k$ is a positive integer, we may replace $\sigma_1$ by $k\sigma_1 = \sigma_1\oplus \sigma_1\oplus \dots \oplus \sigma_1$ ($k$ copies) without changing $r_1$ or the formula to be proved. Likewise for $\sigma_2$. We may therefore assume that $\dim\sigma_1 = \dim\sigma_2$. Next, we choose an indecomposable direct factor $\tau_i$ of $\sigma_i$, $i=1,2$, so as to minimize $\eta(\tau_1\otimes \check \tau_2)$. 
\proclaim{Lemma 1} 
Using the preceding notation, there are positive integers $k$, $a$, $b$ such that 
$$ 
k\sigma_1 = \rho_1\oplus a\tau_1,\quad k\sigma_2 = \rho_2\oplus b\tau_2, 
$$ 
where 
\roster 
\item"\rm (a)" $\dim\rho_1 = \dim\rho_2$ and either 
\item"\rm (b)" $\rho_1$ has no direct factor equivalent to $\tau_1$ or 
\item"\rm (c)" $\rho_2$ has no direct factor equivalent to $\tau_2$. 
\endroster 
\endproclaim 
\demo{Proof} 
Let $m_i$ be the multiplicity of $\tau_i$ in $\sigma_i$ and write $d_i = \dim\tau_i$. By symmetry, we may assume $d_1m_1 \ge d_2m_2$. Thus 
$$ 
d_1\sigma_2 = \rho_2\oplus d_1m_2\tau_2, 
$$ 
for a subspace $\rho_2$ with no factor $\tau_2$. Likewise, 
$$
d_1\sigma_1 = \rho'_1\oplus d_1m_1\tau_1, 
$$
for a subspace $\rho'_1$ with no factor $\tau_1$. We have 
$$ 
\dim\rho_2 - \dim\rho'_1 = d_1^2m_1-d_1d_2m_2. 
$$ 
This integer is divisible by $d_1$ and is non-negative. So, $d_1\sigma_1$ admits a decomposition $d_1\sigma_1 = \rho_1\oplus a\tau_1$ in which $\dim\rho_1 = \dim\rho_2$ and $a$ is a positive integer. The result follows with $k = d_1$ and $b = d_1m_2$. \qed 
\enddemo 
We may replace $(\sigma_1,\sigma_2)$ by $(k\sigma_1,k\sigma_2)$ without changing anything. To simplify notation, we assume that Lemma 1 holds with $k=1$. The hypothesis $r_1r_2 > 1$ implies that one of the spaces $\rho_i$ is non-zero, so both are. Extending notation in the obvious way, we have $r(\rho_1)r(\rho_2) < r_1r_2$ so, by inductive hypothesis, 
$$ 
2\,\eta(\rho_1\otimes\check\rho_2) \ge \eta(\rho_1\otimes\check\rho_1) + \eta(\rho_2\otimes\check\rho_2). 
$$ 
Put 
$$ 
\alpha = \frac{\dim\rho_1}{\dim\sigma_1} = \frac{\dim\rho_2}{\dim\sigma_2},\quad \beta  = 1{-}\alpha. 
$$ 
Applying the definition of $\eta$ to the relations $\sigma_1 = \rho_1\oplus a\tau_1$, $\sigma_2 = \rho_2\oplus b\tau_2$, we get 
$$ 
\align 
\eta(\sigma_1\otimes\check\sigma_2) &= \alpha^2\eta(\rho_1\otimes \check\rho_2) + \alpha\beta\big(\eta(\rho_1\otimes\check\tau_2) {+} \eta(\tau_1\otimes\check\rho_2)\big) + \beta^2\eta(\tau_1\otimes\check\tau_2), \\
\eta(\sigma_1\otimes\check\sigma_1) &= \alpha^2\eta(\rho_1\otimes\check\rho_1) + 2\alpha\beta\eta(\rho_1\otimes\check\tau_1) + \beta^2\eta(\tau_1\otimes\check\tau_1), \\
\eta(\sigma_2\otimes\check\sigma_2) &= \alpha^2\eta(\rho_2\otimes\check\rho_2) + 2\alpha\beta\eta(\rho_2\otimes\check\tau_2) + \beta^2\eta(\tau_2\otimes\check\tau_2). 
\endalign 
$$ 
4.3 Proposition 1 implies that 
$$ 
2\,\eta(\tau_1\otimes\check\tau_2) \ge \eta(\tau_1\otimes\check\tau_1) + \eta(\tau_2\otimes\check\tau_2). 
$$ 
The theorem will therefore follow from: 
\proclaim{Lemma 2} 
With the preceding notation, 
$$ 
\vs(\rho_1\otimes\check \tau_2) + \vs(\tau_1\otimes\check \rho_2) \ge \vs(\rho_1\otimes \check\tau_1) + \vs(\tau_2 \otimes \check \rho_2). 
$$ 
\endproclaim 
\demo{Proof} 
Write $\rho_1 = \bigoplus_{i\in I} \xi_i$ and $\rho_2 = \bigoplus_{j\in J} \theta_j$, where $\xi_i$ and $\theta_j$ are indecomposable. Thus 
$$ 
\eta(\rho_1\otimes\check\tau_1) = \sum_{i\in I} \alpha_i\eta(\xi_i\otimes\check\tau_1), \quad \alpha_i = \dim\xi_i/\dim\rho_1, 
$$ 
and $\sum_{i\in I}\alpha_i = 1$. We have a similar formula for each  of the the three other terms in the inequality to be proved. Combining these, and writing $\beta_j = \dim\theta_j/\dim\rho_2$, the desired relation reduces to 
$$ 
\sum_{i\in I} \alpha_i\eta(\xi_i\otimes\check \tau_2) + \sum_{j\in J}\beta_j\eta(\tau_1\otimes\check \theta_j) 
\ge \sum_{i\in I} \alpha_i\eta(\xi_i\otimes\check \tau_1) + \sum_{j\in J} \beta_j\eta(\tau_2\otimes\check \theta_j). 
$$ 
We multiply each sum over $i$ by $1 = \sum_j\beta_j$ and each in $j$ by $1 = \sum_i\alpha_i$. Comparing the $\alpha_i\beta_j$-term on either side, we see it is enough to prove that 
$$ 
\eta(\xi_i\otimes\check \tau_2) + \eta(\tau_1\otimes\check\theta_j) \ge \eta(\xi_i\otimes\check\tau_1) + \eta(\tau_2\otimes\check \theta_j), \quad i\in I,\ j\in J. 
$$ 
The choice of $(\tau_1,\tau_2)$ gives 
$$ 
\alignedat3 
\eta(\xi_i\otimes\check \tau_2) &\ge \eta(\tau_1\otimes\check\tau_2), &\quad &i\in I, \\
\eta(\tau_1\otimes\check\theta_j) &\ge \eta(\tau_1\otimes\check\tau_2), &\quad &j\in J. 
\endalignedat 
$$
We now apply 4.3 Proposition 2 to get 
$$ 
\alignedat3 
\eta(\xi_i\otimes\check\tau_1) &\le \roman{max}\big\{\eta(\xi_i\otimes\check\tau_2), \eta(\tau_2\otimes\check \tau_1) \big\} &\,&= \eta(\xi_i\otimes\check \tau_2), \\ 
\eta(\tau_2\otimes\check \theta_j) &\le \roman{max}\big\{\eta(\tau_2\otimes\check \tau_1),\eta(\tau_1\otimes\check \theta_j)\big\} &\, &= \eta(\tau_1\otimes\check \theta_j), 
\endalignedat 
$$ 
whence the lemma follows. \qed 
\enddemo 
This completes the proof of 4.4 Theorem and the main assertion of Theorem B. The second assertion of Theorem B is 4.1 Proposition. \qed 
\enddemo 
\subhead 
4.5 
\endsubhead 
To prove Theorem BS, we can pass directly from the end of 4.1 to the start of 4.4. From there on, the argument is identical: one simply replaces $\eta$ by $\vs$ and ``indecomposable'' by ``irreducible'' throughout. 
%%%%%%%%%%%%%%%%%%%%%%%%%%%%%% 
\head\Rm 
5. Upper bounds 
\endhead 
We prove Theorems C and CS. 
\subhead 
5.1 
\endsubhead 
We use a combinatorial device. Let $A = \Bbb Z\big[\Bbb R\big]$ be the integral group ring of the additive group of real numbers. We write the elements of $A$ as finite formal sums of symbols $[\alpha]$, $\alpha\in \Bbb R$. The ring $A$ comes equipped with two canonical homomorphisms 
$$
\aligned 
d: A &\longrightarrow \Bbb Z, \\ [\alpha] &\longmapsto 1, \endaligned \qquad \text{and} \qquad 
\aligned v:A &\longrightarrow \Bbb R, \\ [\alpha] &\longmapsto \alpha. \endaligned 
$$ 
There is a unique bi-additive map $A\times A \to A$, denoted $(x,y) \mapsto x\sve y$, so that 
$$ 
[\alpha]\sve [\beta] = \big[\roman{max}\,\{\alpha,\beta\}\big], \quad \alpha,\beta\in \Bbb R. 
$$ 
Let $A^+$ be the set of elements $\sum_\alpha c_\alpha[\alpha]$ such that $c_\alpha = 0$ if $\alpha<0$ and $c_\alpha \ge 0$ otherwise. 
\proclaim{Proposition} 
If $\sigma,\tau \in A^+$, then 
$$ 
v(\sigma\sve\tau) \le d(\tau)v(\sigma) + d(\sigma)v(\tau) - \roman{min}\,\{v(\sigma),v(\tau)\}. 
$$ 
\endproclaim  
\demo{Proof} 
If either $\sigma$ or $\tau$ is the zero element of $A$, the assertion is trivial. We therefore assume that both $\sigma, \tau \in A^+$ are non-zero and proceed by induction on the integer $d = d(\sigma{+}\tau) \ge 2$. In the first case $d=2$, we have $\sigma = [\alpha]$, $\tau = [\beta]$, for some positive real numbers $\alpha$, $\beta$. The assertion is 
$$ 
\roman{max}\,\{\alpha,\beta\} \le \alpha+\beta-\roman{min}\,\{\alpha,\beta\}. 
$$ 
This holds with equality. For the general inductive step, we may assume by symmetry that $\sigma = \sigma_1{+}\sigma_2$, for non-zero elements $\sigma_i$ of $A^+$. By inductive hypothesis, 
$$ 
v(\sigma_i\sve\tau) \le d(\tau)v(\sigma_i) + d(\sigma_i)v(\tau) - \roman{min}\,\{v(\sigma_i),v(\tau)\}, \quad i=1,2.  
$$ 
Adding and using the inductive hypothesis, we get 
$$ 
\align 
v(\sigma\sve\tau) &\le v(\sigma_1\sve\tau)+v(\sigma_2\sve\tau) \\ 
&\le d(\tau)v(\sigma) + d(\sigma)v(\tau) - \roman{min}\,\{v(\sigma_1),\tau)\} - \roman{min}\,\{v(\sigma_2),\tau)\} \\
&\le d(\tau)v(\sigma) + d(\sigma)v(\tau) - \roman{min}\,\{v(\sigma),v(\tau)\},  
\endalign 
$$ 
as required. \qed 
\enddemo  
\remark{Remark} 
If we fix positive integers $d_1$, $d_2$, and real numbers $v_1$, $v_2$, there exist $\sigma_1,\sigma_2\in A^+$ such that $d_i = d(\sigma_i)$, $v_i = v(\sigma_i)$, and 
$$ 
v(\sigma_1\sve \sigma_2) = d_2v_1+d_1v_2-\roman{min}\,\{v_1,v_2\}. 
$$ 
In other words, the inequality of the proposition is optimal. 
\endremark 
\subhead 
5.2 
\endsubhead 
We prove Theorem CS. Recall that the assertion of the theorem concerning irreducible representations has been proved in 1.3 Lemma. 
\par 
A representation $\sigma\in \wW F$ gives an element $\bold S(\sigma) = \dim(\sigma)[\vs(\sigma)]$ of the ring $A$ of 5.1. For $\sigma\in \wss F$, we define $\bold S(\sigma)\in A^+$ by 
$$ 
\bold S(\sigma_1\oplus \sigma_2\oplus \dots \oplus \sigma_r) = \sum_{i=1}^r \bold S(\sigma_i), \quad \sigma_i \in \wW F. 
$$ 
This definition gives  
$$ 
\aligned 
v\big(\bold S(\sigma)\big) &= \sw(\sigma), \\ 
d\big(\bold S(\sigma)\big) &= \dim\sigma, \endaligned \qquad \sigma \in \wss F. 
$$ 
We know that $\vs(\sigma\otimes\tau) \le \roman{max}\{\vs(\sigma),\vs(\tau)\}$ when both representations $\sigma$, $\tau$ are irreducible. In our present notation, this says 
$$ 
\sw(\sigma\otimes \tau) \le v\big(\bold S(\sigma)\sve \bold S(\tau)\big), \quad \sigma,\tau\in \wW F. 
$$ 
Consequently, if $\rho,\theta\in \wss F$, then 
$$ 
\align 
\sw(\rho\otimes\theta) &\le v\big(\bold S(\rho)\sve \bold S(\theta)\big)  \\ 
&\le d\big(\bold S(\tau)\big)v\big(\bold S(\sigma)\big) + d\big(\bold S(\sigma)\big)v\big(\bold S(\tau)\big) - \roman{min}\,\big \{v\big(\bold S(\sigma)\big),v\big(\bold S(\tau)\big)\big\} \\ 
&= \dim(\tau)\sw(\sigma) + \dim(\sigma)\sw(\tau)-\roman{min}\,\{\sw(\sigma),\sw(\tau)\}, 
\endalign 
$$ 
as required to prove Theorem CS. \qed 
\subhead 
5.3 
\endsubhead 
We can use exactly the same argument to prove Theorem C once we establish: 
\proclaim{Proposition} 
If $\Rho,\Tau \in \dw F$ are indecomposable, then 
$$ 
\eta(\Rho\otimes\Tau) \le \roman{max}\,\{\eta(\Rho),\eta(\Tau)\}. 
$$ 
\endproclaim 
\demo{Proof} 
Let $\Sigma\in \dw F$ be indecomposable. Thus 
$$ 
\eta(\Rho\otimes\Tau) \le \roman{max}\,\{\eta(\Rho\otimes \Sigma), \eta(\check\Sigma\otimes\Tau)\}, 
$$ 
by 4.3 Proposition 2. Taking for $\Sigma$ the trivial character of $\scr W_F$, we get the proposition. \qed 
\enddemo 
This completes the proof of Theorem C. \qed 
%%%%%%%%%%%%%%%%%%%%%%%%%%%%%%%%%%%%%% 
 
\enddocument